\documentclass[11pt,reqno,a4paper]{amsart}

\usepackage{amsmath,amsfonts,amsthm,mathrsfs,amssymb,cite}
\usepackage[usenames]{color} \usepackage{hyperref}

\usepackage{graphicx} \usepackage{latexsym} \usepackage{enumerate}

\theoremstyle{plain} \newtheorem{theorem}{Theorem}[]
\newtheorem*{theorem*}{Theorem} \newtheorem{lemma}[theorem]{Lemma}
\newtheorem*{lemma*}{Lemma} 

\newtheorem*{proposition*}{Proposition}

\newtheorem*{conjecture*}{Conjecture}
\newtheorem{assumption}[theorem]{Assumption}

\theoremstyle{definition} 

\theoremstyle{remark}

\newcommand{\abs}[1]{\left\lvert #1 \right\rvert}
\newcommand{\norm}[1]{\left\lVert #1 \right\rVert}

\newcommand{\R}{\mathbb{R}}

\renewcommand{\eqref}[1]{\textnormal{(\ref{#1})}}

\numberwithin{equation}{section}

\newcommand{\polygon}{\Omega}

\title[Addendum to vanishing transmission eigenfunctions]{Addendum to:
  ``On vanishing near corners of transmission eigenfunctions''}

\author{Eemeli Bl{\aa}sten}

\address{Jockey Club Institute for Advanced Study, Hong Kong
  University of Science and Technology, Hong Kong SAR.}

\email{iaseemeli@ust.hk}

\author{Hongyu Liu} \address{Department of Mathematics, Hong Kong
  Baptist University, Kowloon Tong, Hong Kong SAR.\vspace*{-4mm}}
\address{\vspace*{-4mm}and} \address{HKBU Institute of Research and
  Continuing Education, Virtual University Park, Shenzhen,
  P. R. China.}  \email{hongyu.liuip@gmail.com}

\begin{document}

\maketitle

\begin{abstract}
  In this addendum, we relax a restrictive assumption in
  \cite{BLvanishingCorners} needed for the interior transmission
  eigenfunctions to hold the intrinsic geometric vanishing property in
  a corner. In addition we present in more detail another assumption
  which can also guarantee the vanishing property, namely being
  locally $H^2$ near the corner. This was mentioned briefly in
  \cite{BLvanishingCorners}.
\end{abstract}

\section{Background}
In our recent article on vanishing near corners of transmission
eigenfunctions \cite{BLvanishingCorners}, we proved the following
theorem.

\begin{theorem}[Theorem 3.2 in  \cite{BLvanishingCorners}] \label{thm:1}
  Let $n\in\{2,3\}$ and $V$ be a qualitatively admissible
  potential. Suppose that $k\in\mathbb{R}_+$ is a transmission
  eigenvalue: there exists $v,w \in L^2(\polygon)$ such that
  \begin{align*}
    (\Delta+k^2)v &= 0 \quad \mbox{in}\ \ \polygon,
    \\ (\Delta+k^2(1+V))w &= 0 \quad \mbox{in}\ \ \polygon, \\ w-v
    &\in H^2_0(\polygon), \quad \norm{v}_{L^2(\polygon)} = 1.
  \end{align*}
  If $v$ can be approximated in the $L^2(\polygon)$-norm by a sequence
  of Herglotz waves with uniformly $L^2(\mathbb S^{n-1})$-bounded
  kernels, then
  \begin{equation}\label{eq:01}
  \lim_{r\to +0} \frac{1}{m(B(x_c,r))} \int_{B(x_c,r)} \abs{v(x)} dx =
  0,
  \end{equation}
  where $x_c$ is any vertex of $\polygon$ such that $\varphi(x_c) \neq
  0$.
\end{theorem}

The assumption in Theorem~\ref{thm:1} which says that $v$ can be
approximated by a sequence of Herglotz waves with uniformly bounded
kernels is rather technical. Indeed, since $L^2(\mathbb S^{n-1})$ is a
reflexive Banach space, and exponential functions are in it, then the
assumption above would imply by the weak-* compactness of closed balls
in $L^2(\mathbb S^{n-1})$ that $v$ is actually a Herglotz wave. Hence,
in light of \cite{BPS}, it is a too strong condition. In Section 7 of
\cite{BLvanishingCorners}, we commented on this condition by
describing numerical evidence \cite{numerical} as well as by stating
that using a result in our upcoming paper \cite{BLsingleFarField}, the
theorem can be shown to hold in another alternative general
setting. It is the aim of this addendum to relax the technical
assumption on $v$ in Theorem~\ref{thm:1} to a more reasonable one as
well as to provide the details of the aforementioned general setting for the
theorem as commented in Section 7 of \cite{BLvanishingCorners}.

\section{Relaxed assumptions}

We refer to the assumption in Theorem~\ref{thm:1}, namely that $v$ can
be approximated in the $L^2(\polygon)$-norm by a sequence of Herglotz
waves with uniformly $L^2(\mathbb S^{n-1})$-bounded kernels, as
Assumption~1. Next, we provide two alternative assumptions, and either
of them in place of Assumption~1 leads to the same conclusion in
Theorem~\ref{thm:1}.

\begin{assumption} \label{ass2}
  The transmission eigenfunction $v$ can be approximated in the
  $L^2(\Omega)$-norm by Herglotz waves $v_j$, $j\in{({1,\infty})}$
  with kernels $g_j$ such that
  \begin{equation}
    \norm{v-v_j}_{L^2(\Omega)} < e^{-j}, \qquad
    \norm{g_j}_{L^2(\mathbb S^{n-1})} \leq C (\ln j)^\beta
  \end{equation}
  for some positive constants $C$ and $0<\beta<1/(2n+8)$,
\end{assumption}

\begin{assumption} \label{ass3}
  The transmission eigenfunctions $v,w$ are $H^2$-smooth in a
  neighbourhood of $x_c$ in $\Omega$.
\end{assumption}

Clearly, Assumption~\ref{ass2} is much weaker than Assumption~1, and
it allows the kernels of the Herglotz waves to blow up under a certain
logarithmic rate. Assumption~\ref{ass3} is simply a condition on the
boundary: elliptic regularity guarantees $H^2$-smoothness in any
domain in $\Omega$ whose boundary is disjoint from $\partial\Omega$.
Assumption~\ref{ass2} is generally true based on numerical evidence
\cite{numerical}. That paper also shows that Assumption~\ref{ass3}
does not always hold: if the corner is not convex then actually $v$
and $w$ blow up. This observation screams for a mathematical proof.

\medskip
Let us first show that Assumption~\ref{ass3} implies the conclusion of
our theorem. The following auxiliary result shall be needed.

\begin{lemma}[Lemma~3.5 in \cite{BLsingleFarField}] \label{boundary2corner}
  Let $\tilde\Omega\subset\R^n$, $n\in\{2,3\}$ be a bounded
  domain. Let $q,q'\in L^\infty(\tilde\Omega)$ and let $u, u' \in
  H^2(\tilde\Omega)$ solve
  \[
  (\Delta+q)u = 0, \qquad (\Delta+q')u' = 0
  \]
  in $\tilde\Omega$. Let $x_c \in \partial\tilde\Omega$ be a point for
  which there is a neighbourhood $B$ and an admissible cell $\Sigma$
  having $x_c$ as vertex such that $B\cap\tilde\Omega=B\cap\Sigma$.

  Assume that $u=u'$ and $\partial_\nu u = \partial_\nu u'$ on
  $B\cap\partial\tilde\Omega$. If $q$ and $q'$ are $C^\alpha$
  uniformly H\"older-continuous functions with $\alpha>0$ in 2D and
  $\alpha>1/4$ in 3D in $B \cap \tilde\Omega$ then
  \[
  (q-q')(x_c) u(x_c) = (q-q')(x_c) u'(x_c) =0.
  \]
\end{lemma}

Using Lemma~\ref{boundary2corner}, the proof of Theorem~\ref{thm:1}
with Assumption~1 replaced by Assumption~\ref{ass3} can proceed as
follows. Set $q = k^2(1+V)$, $q'=k^2$ and $u=w$, $u'=v$. Moreover set
$\tilde\Omega$ as the neighbourhood of $x_c$ of Assumption~\ref{ass3}
where $v,w$ are $H^2$-smooth. Since $v-w\in H^2_0(\Omega)$ then $v=w$
and $\partial_\nu v = \partial_\nu w$ on a suitable part of the
boundary of $\tilde\Omega$, and the lemma readily implies
$v(x_c)=w(x_c)=0$. That is, the vanishing property holds at the corner
point.

\medskip

Next, we consider Theorem~\ref{thm:1} with Assumption~1 replaced by
Assumption~\ref{ass2}. We shall make modifications to its proof,
i.e. to \emph{Proof~of~Theorem~3.2} in Section 6 of
\cite{BLvanishingCorners} as follows. Let $v_j$ be the sequence of
Herglotz waves given by Assumption~\ref{ass2}. If there is a
subsequence, again denoted by $v_j$, with orders $N_j\geq1$,
i.e. $v_j(x_c)=0$, then the Taylor expansion gives
\[
\norm{v_j}_{L^1(B(x_c,r))} \leq 0 + C_{k,n} \norm{g_j}_{L^2(\mathbb
  S^{n-1})} r^{n+1}.
\]
Using $\lVert{v-v_j}\rVert_{L^1(B(x_c,r)\cap\Omega)} \leq C_\Omega
e^{-j}$, one has
\[
\frac{1}{m(B(x_c,r)\cap\Omega)} \int_{B(x_c,r)\cap\Omega} \abs{v(x)}
dx \leq C_{k,n,\Omega} ( e^{-j} r^{-n} + \norm{g_j}_{L^2(\mathbb
  S^{n-1})} r ).
\]
Since $\lVert{g_j}\rVert_{L^2(\mathbb S^{n-1})} \leq C (\ln j)^\beta$,
we may choose $j = (n+1) \ln \frac{1}{r}$ to see that the average of
$\lvert{v}\rvert$ over $B(x_c,r)\cap\Omega$ tends to zero as $r\to
+0$.

\smallskip
For the other case, we may assume that $v_j(x_c)\neq0$ for all $j$, in
other words $N_j=0$. Then $v_j(x) = v_j(x_c) + r_j(x)$ with
$\lvert{r_j(x)}\rvert \leq C_{k,n} \lVert{g_j}\rVert_{L^2(\mathbb
  S^{n-1})} \lvert{x-x_c}\rvert$ as above. The trivial upper bound
$\lvert{v_j(x_c)}\rvert \leq C \lVert{g_j}\rVert_{L^2(\mathbb
  S^{n-1})}$ is not strong enough to conclude the vanishing of
$v(x_c)$. The main motivation and idea of our paper
\cite{BLvanishingCorners} start next.

Using Theorem~3.1 of \cite{BLvanishingCorners} with the incident wave
$v_j/\lVert{g_j}\rVert_{L^2(\mathbb S^{n-1})}$, one can obtain the
following lower bound of the corresponding far-field pattern
\begin{equation}\label{eq:1}
\norm{v^s_{j \infty}} \geq \frac{\mathcal{S} \norm{g_j}_{L^2(\mathbb
    S^{n-1})}}{\exp \exp c \min (1,
  \frac{\norm{P_{j,N_j}}}{\norm{g_j}_{L^2(\mathbb S^{n-1})}})^{-\ell}}
\end{equation}
where $\mathcal S = \mathcal S(V,k)$, $c = c(V,k,n)$ and $\ell =
2(n+4)$. The theorem statement is unclear about $\ell$ but its proof
allows this choice. Moreover in this case $N_j=0$ for all $j$ so
$P_{j,N_j} = v(x_c)$ and thus $\lVert{P_{j,N_j}}\rVert = c_n
\lvert{v_j(x_c)}\rvert$. The critical observation for the transmission
eigenfunctions to vanish at corners is to note that since
$v^s_{j\infty}\to0$ (Proposition~4.2 in \cite{BLvanishingCorners})
then also $v_j(x_c)\to0$, and so \eqref{eq:01}, as follows.

One can first use the estimate $\lVert{g_j}\rVert_{L^2(\mathbb
  S^{n-1})} \geq c_{n,\Omega} > 0$ for $j$ large enough in the
numerator of \eqref{eq:1}. Indeed, this is implied by
\[
\norm{v-v_j}_{L^2(\Omega)}<e^{-j}, \quad \norm{v}_{L^2(\Omega)}=1,
\quad \norm{v_j}_{L^2(\Omega)} \leq C_{\Omega,n}
\norm{g_j}_{L^2(\mathbb S^{n-1})}.
\]
Next, by Proposition~4.2 in \cite{BLvanishingCorners}, one can
estimate $\lVert{v_{j\infty}^s}\rVert_{L^2(\mathbb S^{n-1})} < C_{V,k}
e^{-j}$ in the left hand side of \eqref{eq:1}. After our estimates and
recalling that $\lVert{P_{j,N_j}}\rVert = c_n\lvert{v_j(x_c)}\rvert$,
Equation \eqref{eq:1} becomes
\[
\exp\left( -j + \exp c \min \left(1, \frac{c_n
  \abs{v_j(x_c)}}{\norm{g_j}_{L^2(\mathbb S^{n-1})}} \right)^{-\ell}
\right) \geq C_{V,k}^{-1} \mathcal{S} c_{n,\Omega}.
\]
Let us write $j_0 = \ln(C_{V,k}^{-1} \mathcal{S} c_{n,\Omega})$ and $z
= \min (1, c_n \lvert{v_j(x_c)}\rvert / \lVert{g_j}\rVert_{L^2(\mathbb
  S^{n-1})})$. Now, when $j+j_0>0$, we have $z^{-\ell} \geq
c^{-1}\ln(j+j_0)$ and consequently
\[
z \leq c^{1/\ell} \big(\ln(j+j_0)\big)^{-1/\ell}.
\]
If $j$ is large enough then the right-hand side is smaller than
$1$. Hence the minimum in the definition of $z$ gives
\[
\abs{v_j(x_c)} \leq (c/c_n)^{1/\ell} \norm{g_j}_{L^2(\mathbb S^{n-1})}
\big(\ln (j+j_0)\big)^{-1/\ell}
\]
for large $j$. Finally, writing $B=B(x_c,r)$ and recalling $v_j(x) =
v_j(x_c) + r_j(x)$, one has
\begin{equation}\label{eq:2}
\begin{split}
  &\frac{\norm{v}_{L^1(B\cap\Omega)}}{m(B\cap\Omega)} \leq C_n r^{-n}
  \big( \norm{v-v_j}_{L^1(B\cap\Omega)} + \norm{v_j(x_c) +
    r_j(x)}_{L^1(B\cap\Omega)} \big) \\&\qquad \leq C \big( e^{-j}
  r^{-n} + \norm{g_j}_{L^2(\mathbb S^{n-1})} (\ln (j+j_0))^{-1/\ell} +
  \norm{g_j}_{L^2(\mathbb S^{n-1})} r \big) \\ &\qquad \leq C
  \big(e^{-j}r^{-n} + (\ln (j+j_0))^{\beta-1/\ell} + (\ln j)^\beta r
  \big).
\end{split}
\end{equation}
By choosing $j=(n+1) \ln \frac{1}{r}$ with $r$ small enough and noting
that $\beta < 1/\ell = 1/(2n+8)$, one readily sees that the right-hand
side of \eqref{eq:2} tends to zero as $r\to +0$, which in turn implies
\eqref{eq:01}.

\end{document}